\documentclass{article}

\usepackage{arxiv}

\usepackage{amsmath,amsthm,amssymb,amscd}
\usepackage[all,cmtip]{xy}
\usepackage{diagxy}
\usepackage[utf8]{inputenc} 
\usepackage[T1]{fontenc}    
\usepackage{hyperref}       
\usepackage{url}            
\usepackage{booktabs}       
\usepackage{amsfonts}       
\usepackage{nicefrac}       
\usepackage{microtype}      
\usepackage{lipsum}
\usepackage{mathrsfs}
\usepackage{fancyhdr}
\usepackage[dvips]{color}
\usepackage{color}
\usepackage{subfig}
\usepackage{graphicx}

\newtheorem{exam.}[subsection]{Example}
\newtheorem{def.}[subsection]{Definition}
\newtheorem{Pdef}[subsection]{Proposition/Definition}

\newtheorem{lemma}[subsection]{Lemma}
\newtheorem{prop.}[subsection]{Proposition}
\newtheorem{coro}[subsection]{Corollary}
\newtheorem{theorem}[subsection]{Theorem}

\newtheorem*{remark}{Remark}

\newcommand{\id}{{\rm{id}}}

\newcommand{\bbZ}{\mathbb{Z}}
\newcommand{\bbR}{\mathbb{R}}
\newcommand{\bbN}{\mathbb{N}}
\newcommand{\bbC}{\mathbb{C}}

\newcommand{\shO}{\mathscr{O}}

\newcommand{\shF}{\mathscr{F}}

\newcommand{\shD}{\mathscr{D}}
\newcommand{\shG}{\mathscr{G}}

\newcommand{\shC}{\mathscr{C}}

\title{A Brief Note for Sheaf Structures on Posets}
\author{
  Chuan-Shen Hu \\ 
  Department of Mathematics\\
  National Taiwan Normal University\\
  Taipei, Taiwan \\
  \texttt{peterbill26@hotmail.com}
}

\begin{document}
\maketitle

\begin{abstract}
This note is a part of the lecture notes of a graduate student algebraic geometry seminar held at the department of mathematics in National Taiwan Normal University, 2020 Falls. It aims to introduce an example of sheaves defined on posets equipped with the \textit{Alexandrov topology}, called the \textit{cellular sheaves}. A cellular sheaf is a functor from the category of a poset to the category of specific algebraic structures (\textit{e.g.} the category of groups). Strictly speaking, even equipping the poset with the Alexandrov topology, it is just the definition of a \textit{pre-sheaf} on the Alexandrov topological space. By checking details, cellular sheaves are actually sheaves on topological spaces. This is a well-known fact in sheaf theory via the Kan extension, while it requires readers who are familiar with the category theory. In this note, we follow an elementary approach to describe the connection between cellular sheaves and sheaves concerned in algebraic geometry, where only basic commutative algebra and point-set topology are required as the background knowledge.   
\end{abstract}

\section{Introduction}
\label{Intro}
Sheaves on topological spaces are the main targets in modern algebraic geometry~\cite{Hartshorne, Ueno, LeiFu, QingLiu}, which assigns each open subset of a topological space to an algebraic structure and connects them via homomorphisms. In the note, all rings are assumed to be commutative with identity elements, and all homomorphisms of rings are assumed to be identity-preserving.

Let $X$ be a topological space, a \textit{pre-sheaf} $\shF$ of Abelian groups (resp. rings, algebras, vector spaces, ... etc.) on $X$(cf.~\cite{bredon1997sheaf, Hartshorne}) consists the following data:
\begin{itemize}
\item For every open subset $U$ of $X$, there is an Abelian group (resp. ring, algebra, vector space, ... etc.) $\shF(U)$, in particular, $\shF(\emptyset) = \lbrace 0 \rbrace$. Elements in $\shF(U)$ are called \textit{sections} of $\shF$ on $U$.
\item For every inclusion $V \subseteq U$ of open subsets of $X$, there is a group (resp. ring, algebra, vector space, ... etc.) homomorphism $\rho_{UV} : \shF(U) \rightarrow \shF(V)$, called the \textit{restriction map}.
\end{itemize}
The data satisfy the following two properties :
\begin{itemize}
\item[(1)] For any open subset $U$ of $X$, $\rho_{UU} = \id_{\shF(U)}$.
\item[(2)] For every inclusions $W \subseteq V \subseteq U$ of open subsets of $X$, $\rho_{UW} = \rho_{VW} \circ \rho_{UV}$. In other words, for any inclusions $W \subseteq V \subseteq U$, the following diagram is commutative.
\begin{equation}
\label{Eq : Diagram of Restriction of maps of Sheaves}
\xymatrix{
\shF(U) \ar[r]^{\rho_{UV}} \ar[dr]_{\rho_{UW}} & \shF(V) \ar[d]^{\rho_{VW}} \\
& \shF(W)
}
\end{equation}  
\end{itemize}
Intuitively, the purpose for defining sections $s \in \shF(U)$ is trying to copy the concept of functions on open set $U$ and view $\shF(U)$ as the generalized space of functions on $U$. Hence, for section $s \in \shF(U)$, we often denote $\rho_{UV}(s)$ by $s|_V$. A presheaves $\shF$ of Abelian groups on $X$ is a \textit{sheaves} if it satisfies the following conditions :
\begin{itemize}
\itemsep = -1pt
\item[(a)] Let $s, t \in \shF(U)$ be sections. If $\lbrace U_i \rbrace_{i \in I}$ is an open covering of $U$ (i.e., $U = \bigcup_{i \in I} U_i$) such that $s|_{U_i} = t|_{U_i}$ for all $i \in I$, then $s = t$.
\item[(b)] Suppose $\lbrace U_i \rbrace_{i \in I}$ is an open covering of $U$ and $s_i \in \shF(U_i)$ are some sections such that $s_i|_{U_i \cap U_j} = s_j|_{U_i \cap U_j}$ whenever $i, j \in I$. Then there exists a section $s \in \shF(U)$ such that $s|_{U_i} = s_i$ for any $i \in I$. We note that condition (a)  implies that the existence of such $s$ is unique. 
\end{itemize}
For presheaves of Abelian groups, we can use the language of exact sequences of Abelian groups to describe the definition of sheaves.
\begin{prop.}
[\cite{QingLiu}~Lemma 2.7]
A presheaf of Abelian groups on $X$ is a sheaf if and only if for any open covering $\lbrace U_i \rbrace_{i \in I}$ of any open subsets $U$, the sequence 
\begin{equation}
\label{Equation : Exact Form for Presheaf Being a Sheaf}
0 \longrightarrow \shF(U) \stackrel{\varphi}\longrightarrow \prod_{i \in I} \shF(U_i) \stackrel{\psi}\longrightarrow \prod_{(i,j) \in I^2} \shF(U_i \cap U_j)
\end{equation}
is exact, where $\varphi$ represents the mapping $ s \longmapsto (s|_{U_i})_{i \in I}  $ 
and $\psi$ is the mapping $(s_i)_{i \in I} \longmapsto (s_j|_{U_i \cap U_j} - s_i|_{U_i \cap U_j})_{(i,j) \in I^2}$.
\end{prop.}
\begin{exam.}
Let $X$ be a topological space. For every non-empty open subset $U$ of $X$, set $\shC(X)$ to be the set of all continuous functions from $U$ to $\bbC$ and $\shC(\emptyset) = \lbrace 0 \rbrace$. On the other hand, if $V \subseteq U$ are non-empty open subsets of $X$, then $\rho_{UV}(f)$ is defined to be the restriction function $f|_V$ whenever $ f \in \shC(U)$. Then $\shC$ is a sheaf of rings on $X$.
\end{exam.}
In classical algebraic geometry, sheaf is a natural structure for describing regular functions defined on Zariski open subsets of a variety over an algebraically closed field. 
\begin{exam.}
[\cite{Hartshorne}~Example II.1.0.1]
Let $k = \overline{k}$ be an algebraically closed field and $X$ be an affine variety over $k$. For each open set $U$ of $X$, we define $\shO_X(U)$ to be the $k$-algebra of all regular functions on $U$. By equipping restriction of functions among all open sets of $X$, $\shO_X$ is a sheaf of $k$-algebras on $X$.
\end{exam.}
Here we see a well-known example about a presheaf which is not a sheaf.
\begin{exam.}
Let $X = \bbR$ be equipped with the usual topology. Define $\shF$ on $X$ by $\shF(U) = \{ f : U \rightarrow \bbR : f \ {\rm is \ bounded} \}$ equipped with function restrictions. Then $\shF$ is a presheaf of Abelian groups but $\shF$ is not a sheaf.
\end{exam.}
\begin{proof}
It is clear that $\shF$ is a presheaf of Abelian groups since sum of bounded functions is still bounded. Moreover, it also satisfies (a). To prove $\shF$ doesn't satisfy (b), we consider open intervals $U_n = (-n,n)$, $n \in \bbN$. Consider $f : X \rightarrow \bbR$ by $f(x) = x$ and $f_n = f|_{U_n}$. Clearly, $\bbR = \bigcup_{n = 1}^{\infty} U_n$, $f_n \in \shF(U_n)$ and $f_n|_{U_n \cap U_m} = f_m|_{U_n \cap U_m}$ whenever $m, n \in \bbN$. However, the gluing function of all $f_n$ must be $f$ but $f$ is not bounded.
\end{proof} 
On the other hand, a \textit{cellular sheaf} $\shG$~\cite{ghrist2014elementary, robinson2014topological, RobinsonHunting, Curry2015, justinCurryPHDThesis} defined on a poset (i.e. a partially ordered set) $(P, \leq)$ is a rule consists of the following data:
\begin{itemize}
    \item For every $p \in P$, there is an Abelian group (resp. rings, algebras, vector spaces, ... etc.) $\shG(p)$, in particular, $\shG(\emptyset) = \lbrace 0 \rbrace$. Elements in $\shG(p)$ are called \textit{sections} of $\shG$ on $U$.
    \item For every inclusion $p \leq q$ in $P$, there is a group (resp. rings, algebras, vector spaces, ... etc.) homomorphism $\rho_{pq} : \shG(p) \rightarrow \shG(q)$, called the \textit{restriction map}. 
\end{itemize}
The data satisfy the following two properties :
\begin{itemize}
\item[(1)] For any open subset $U$ of $X$, $\rho_{pp} = \id_{\shG(p)}$.
\item[(2)] For every $p \leq q \leq r$ in $P$, $\rho_{pr} = \rho_{qr} \circ \rho_{pq}$. In other words, the diagram
\begin{equation}
\label{Eq : Diagram of Restriction of maps of Cellular Sheaves}
\xymatrix{
\shG(p) \ar[r]^{\rho_{pq}} \ar[dr]_{\rho_{pr}} & \shG(q) \ar[d]^{\rho_{qr}} \\
& \shG(r)
}
\end{equation} 
is commutative for every $p \leq q \leq r$.
\end{itemize}
Compared with sheaves defined on topological spaces, a pre-sheaf $\shF$  over a topological space $X$ can be viewed as a cellular sheaf from the poset $({\rm Open}(X), \leq)$ to the category of Abelian groups where $U \leq V$ if and only if $V \subseteq U$. Conversely, given a cellular sheaf $\shG$ on a poset $(P, \leq)$, we may also construct a topology $\mathfrak{A}$, called the \textit{Alexandrov topology} on $P$. Then $\shG$ can be canonically encoded as a pre-sheaf (still denoted by $\shG$) of Abelian groups (resp. rings, modules ... etc.) over the topological space $(P,\mathfrak{A})$. Actually, the extended pre-sheaf $\shG$ on $(P,\mathfrak{A})$ is a sheaf. This fact can be proved via categorical languages~\cite{justinCurryPHDThesis, Ladkani_2008}. Concrete construction for the extension were also given in~\cite{ghrist2020cellular, RobinsonHunting, robinson2014topological, RobinsonNyquist, Curry2015}. Equivalent to the concrete approach, we use the base sheaf (i.e., $\mathfrak{B}$-sheaf) structure to express the extension, which is much more familiar with people who are using the language in algebraic geometry~\cite{Hartshorne, QingLiu, AM}. The same approach about using $\mathfrak{B}$-sheaves can be also found in~\cite{RonHeld1}, while the requirements of background knowledge in category theory (\textit{e.g.} Abelian categories, limits) are still essential, and details for $\mathfrak{B}$-sheaf structures were left as exercises. In this note, we fulfill all elementary and essential details about Alexandrov topology and $\mathfrak{B}$-sheaf structure. On the other hand, we focus on the category of $A$-modules as the target category. It is naturally an Abelian category and the (co)limits exist in the category, hence many categorical preliminaries could be avoided. 
\section{Morphisms and Stalks}
Let $(X,\mathfrak{T})$ be a topological space and $\shF$ a presheaf of Abelian groups (resp. rings, ... etc.)  on $X$, where $\mathfrak{T}$ is the topology of all open subsets of $X$. Then the relation $\leq$ defined on $\mathfrak{T}$ by
\begin{equation*}
U \leq V \ \ {\rm if \ and \ only \ if } \ \ V \subseteq U
\end{equation*} 
is a partial order on $\mathfrak{T}$. Moreover, the pair $(\mathfrak{T}, \leq)$ forms a direct set since $ U, V \leq U \cap V $ whenever $ U, V \in \mathfrak{T}$. Especially, if $P \in X$ is a points, then the pair $(\mathfrak{T}_P, \leq)$ also forms a direct set, where $\mathfrak{T}_P := \lbrace U \in \mathfrak{T} \ | \ P \in U \rbrace$. We define the \textit{stalk} $\shF_P$ of $\shF$ at $P$ to be the direct limit~\cite{AM}
\begin{equation*}
\shF_P = \left. \left( \bigoplus_{U \in \mathfrak{T}_P} \shF(U) \right)  \right/ \shD_P,
\end{equation*}
where $\shD_P$ is the submodule generated by all elements of the form $s - \rho_{UV}(s)$, where $s \in \shF(U)$ and $U, V \in \mathfrak{T}_P$ with $ V \subseteq U$ (i.e., $U \leq V$). For $s \in \shF(U)$, we often use the notation $s_P = s + \shD_P$ to denote the image of the canonical map $\shF(U) \rightarrow \shF_P$, which is called the \textit{germ} of $s$ at $P$. By the structure of direct limits, every element in  $\shF_P$ can be represented by $s_P$, where $s$ is a section in $\shF(U)$ with $U \in \mathfrak{T}_P$. Moreover, if $\shF$ is a presheaf of rings on $X$, then $\shF_P$ is also ring, where the multiplication is defined by
\begin{equation*}
s_P \cdot t_P = (s|_{U \cap V} \cdot t|_{U \cap V})_P
\end{equation*}
whenever $s \in \shF(U)$, $t \in \shF(V)$ with $U, V \in \mathfrak{T}_P$. For a presheaf $\shF$ of rings on $X$, the stalk $\shF_P$ of $\shF$ at $P$ is defined by its group structures on each $\shF(U)$ with $U \in \mathfrak{T}_P$. Therefore, in many situations, to investigate the germs of sections of a sheaf $\shF$ of rings (or other objects with based group structure), it is sufficient to investigate the germs of sections of the sheaf of Abelian groups inherits from $\shF$.
\begin{prop.}
\label{Proposition : Stalks on section determine the section}
Let $X$ be a topological space and $\shF$ be a sheaf of Abelian groups on $X$. If $U$ is an open subset of $X$, $s, t \in \shF(U)$ and $s_P = t_P$ for all $P \in U$, then $s = t$.
\end{prop.}
\begin{proof}
If $(s - t)_P = 0$ in $\shF_P$, then there is an open neighbourhood $U_P \subseteq U$ of $P$ in $X$ such that $(s-t)|_{U_P} = 0$ in $\shF(U_P)$. Because $\lbrace U_P \rbrace$ is an open covering of $U$ and $\shF$ is a sheaf, (a) shows that $s - t = 0$ in $\shF(U)$, as desired.
\end{proof}
Given any two algebraic objects, it is natural to consider homomorphisms between them. In sheaf theory, the analogous concepts are known as the morphisms of sheaves.

\begin{def.}
Let $\shF$ and $\shG$ be presheaves of Abelian groups (resp. rings, algebras, vector spaces, ... etc.) on $X$. A \textbf{morphism of presheaves} $\phi : \shF \rightarrow \shG$ consists of groups (resp. rings, algebras, vector spaces, ... etc.) homomorphism $\phi(U) : \shF(U) \rightarrow \shG(U)$ for every open subset $U$ such that for every $V \subseteq U$ of open subsets, the diagram
\begin{equation}
\label{Equation : Commutativity Diagram of Morphisms}
\bfig \Square(2000,0)|aaaa|[\shF(U)`\shG(U)`\shF(V)`\shG(V);\phi(U)`\rho_{UV}`\sigma_{UV}`\phi(V)]\efig 
\end{equation}
commutes. For any $P \in X$, $\phi$ induces a homomorphism on stalks $\phi_P : \shF_P \rightarrow \shG_P$ by $s_P \mapsto \phi(U)(s)_P$ whenever $s \in \shF(U)$. For any presheaf $\shF$, we have the identity morphism $\id_{\shF}$ i.e., $\id_{\shF}(U) = \id_{\shF(U)}$ for each $U \in \mathfrak{T}(X)$.
\end{def.}
\begin{def.}
A morphism of presheaves $\phi : \shF \rightarrow \shG$ is called an \textbf{isomorphism} if it has two-sided inverses, that is, there is a morphism of presheaves $\psi : \shG \rightarrow \shF$ such that $\psi \circ \phi = \id_{\shF}$ and $\phi \circ \psi = \id_{\shG}$. Obviously, it is equivalent to say that $\phi(U) : \shF(U) \rightarrow \shG(U)$ is an isomorphism for every open subset $U$ of $X$.
\end{def.}
\begin{theorem}
\label{Theorem : phi of sheaves is an iso iff every phi(U) is an iso}
Let $\phi : \shF \rightarrow \shG$ be a morphism of sheaves on a topological space $X$. Then $\phi$ is an isomorphism if and only if the induced map on stalks $\phi_P : \shF_P \rightarrow \shG_P$ is an isomorphism for every $P \in X$.
\end{theorem}
\begin{proof}
See Proposition 1.1 in~\cite{Hartshorne} or Proposition 1.1.1 in~\cite{LeiFu}.
\end{proof}
\begin{def.}
Let $X$ be a topological space and $\shF, \shG$ be sheaves of Abelian groups. A morphism $\phi : \shF \rightarrow \shG$ is called \textbf{injective} (resp. \textbf{surjective}) if $\phi_P : \shF_P \rightarrow \shG_P$ is injective (resp. surjective) for all $P \in X$.
\end{def.}
\begin{coro}
Let $X$ be a topological space and $\shF, \shG$ be sheaves of Abelian groups. A morphism $\phi : \shF \rightarrow \shG$ is injective if and only if $\phi(U) : \shF(U) \rightarrow \shG(U)$ is injective for every open subset $U$ of $X$.
\end{coro}
\begin{proof}
The direction $\Rightarrow$ was proved in Theorem \ref{Theorem : phi of sheaves is an iso iff every phi(U) is an iso}. Conversely, if $\phi(U)$ is injective for every $U$ and $s_P \in \ker(\phi_P)$ with $s \in \shF(U)$ and $P \in U$, then $0 = \phi_P(s_P) = \phi(U)(s)_P$. This shows that $\phi(V)(s|_V) = 0$ for some open neighbourhood $V$ of $P$ in $U$. Because $\phi(V)$ is one-to-one, $s|_V = 0$ and we conclude $s_P = 0$.
\end{proof}

\section{Posets and Alexandrov Topology}
\label{sec:Posets and Alexandrov Topology}
To introduce the Alexandrov Topology, we first recall the definition of pre-ordered set and partially ordered set. Further information can be found in~\cite{Ale37, Ale47}.
\begin{def.}
A non-empty set $P$ which is equipped with a relation $\leq$ on $P$ is called a \textbf{pre-ordered set} if the relation $\leq$ satisfies the following two properties:
\begin{itemize}
    \item (Reflexive) $x \leq x$ whenever $x \in P$;
    \item (Transitive) If $x \leq y$ and $y \leq z$, then $x \leq z$.
\end{itemize}
In addition, a pre-ordered set $(P,\leq)$ is called a \textbf{partially ordered set} (or, just a \textbf{poset}) if $\leq$ satisfies the property:
\begin{itemize}
    \item (Anti-symmetric) For $x, y \in P$, if $x \leq y$ and $y \leq x$, then $x = y$.
\end{itemize}
\end{def.}
It is clear that $(\bbZ, \leq)$, $(\bbR, \leq)$, ... etc. are posets. Except posets of number systems, the following are some well-known examples we used frequently.
\begin{exam.}
Let $X$ be a non-empty set and $2^X$ be its power set, then the set inclusion relation $\subseteq$ on $2^X$ makes $(2^X,\subseteq)$ be a poset.
\end{exam.}
\begin{exam.}
Let $(P,\leq)$ be a poset and $P^{\bbN}$ be the set of sequences $(a_n)_{n \in \bbN}$ in $P$. Then the relation $\leq$ defined by $(a_n)_{n \in \bbN} \leq (b_n)_{n \in \bbN}$ if and only if $a_n \leq b_n$ for all $n \in \bbN$ is a partial order on $P^{\bbN}$.
\end{exam.}  
\begin{exam.}
Let $X$ be an arbitrary set and $(P,\leq)$ be a poset. Suppose $P^X$ is the set of all functions from $X$ to $P$, then the relation $\leq$ on $P^X$ defined by $f \leq g$ if and only if $f(x) \leq g(x)$ for all $x \in X$ is a partial order.  
\end{exam.}
\begin{exam.}
[Abstract simplicial complex~\cite{munkres2018elements}]
Let $V$ be a non-empty set. A subset $\mathcal{K}_V$ of $2^V$ is called an \textbf{abstract simplicial complex with vertices set} $V$ if it satisfies the following properties:
\begin{itemize}
\item If $\sigma \in \mathcal{K}_V$, then $\sigma \neq \emptyset$ and $\# \sigma < \infty$;
\item If $\sigma \in \mathcal{K}_V$ and $\emptyset \neq \tau \subseteq \sigma$, then $\tau \in \mathcal{K}_V$.
\end{itemize}
Each $\sigma$ is called a $\textbf q$-\textbf{simplex} where $q = \#\sigma - 1$. We can define relation $\preceq$ on $\mathcal{K}_V$ by $\tau \preceq \sigma$, called $\tau$ is a \textbf{face} of $\sigma$ if and only if $\tau \subseteq \sigma$. By inheriting the order property of $\subseteq$ on $2^V$, the pair $(\mathcal{K}_V,\preceq)$ is a poset.
\end{exam.}
We also exhibit an example of $\leq$ on a set $X$ such that $(X,\leq)$ is a pre-ordered set but it is not a poset.
\begin{exam.}
Let $\bbN$ be the set of all positive integers and $\mathcal{P} \subseteq 2^\bbN$ be the set of all finite subsets of $\bbN$. Define relation $\preceq$ on $\mathcal{P}$ by $A \preceq B$ if and only if $|A| \leq |B|$. Then $(\mathcal{P},\preceq)$ is a pre-ordered set but it is not a poset.
\end{exam.}
\begin{proof}
It is clear that $A \preceq A$ for every $A \in \mathcal{P}$. Moreover, because $(\bbZ,\leq)$ is a pre-ordered set, we have $|A| \leq |B|$ and $|B| \leq |C|$ implies that $|A| \leq |C|$, this shows that $A \preceq B$ and $B \preceq C$ implies $A \preceq C$. Hence $(\mathcal{P},\preceq)$ is a pre-ordered set. However, if $A = \{ 1,2\}$ and $B = \{ 2,3\}$, then  $A,B \in \mathcal{P}$ and $A \preceq B$, $B \preceq A$ since $|A| = |B| = 2$. However, because $A \neq B$, we deduce that $(\mathcal{P},\preceq)$ is not a poset. 
\end{proof}
\begin{def.}
Let $(P,\leq)$ and $(Q,\preceq)$ be pre-ordered sets. A function $f : P \rightarrow Q$ is called \textbf{order-preserving} if $f(x) \preceq f(y)$ whenever $x \leq y$.
\end{def.}
\begin{prop.}
[Poset associated by a pre-ordered set]
Let $(P,\leq)$ be a pre-ordered set, then the relation $\sim$ on $P$ defined by $x \sim y$ if and only if $x \leq y$ and $y \leq x$ is an equivalence relation. If $\overline{P}$ is the set of equivalence classes with respect to $\sim$, then the relation $\leq$ on $\overline{P}$ defined by $\overline{x} \leq \overline{y}$ if and only if $x\leq y$ is a well-defined partial order on $\overline{P}$ and the canonical mapping $\pi : P \rightarrow \overline{P}$ is order-preserving. 

Moreover, we have the following universal property: If $(Q,\preceq)$ is a pre-ordered set and $f : (P,\leq) \rightarrow (Q,\preceq)$ is order-preserving, then there is a unique order-preserving map $\overline{f} : (\overline{P},\leq) \rightarrow (Q,\preceq)$ such that $\overline{f} \circ \pi = f$.
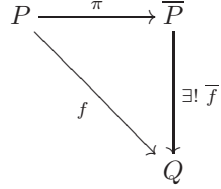
\begin{equation*}
\xymatrix@+2.0em{
				& P
				\ar[dr]_{f}
				\ar[r]^{\pi}
                & \overline{P} 
                \ar[d]^{\exists! \ \overline{f}}
                \\
        		& 
				& Q
}
\end{equation*}
\end{prop.}
\begin{proof}
We first check that $\sim$ is an equivalence relation on $P$. Because $x \leq x$ for every $x \in P$, we have $x \sim x$. If $x \sim y$, then both $x \leq y$ and $y \leq x$. By definition of $\sim$, we have $y \sim x$. Finally, if $x \sim y$ and $y \sim z$, then $x \leq y$, $y \leq x$, $y \leq z$ and $z \leq y$. By the transitivity of $\leq$, we deduce that $x \leq z$ and $z \leq x$.

If $\overline{x} = \overline{x'}$ and $\overline{y} = \overline{y'}$, then $x \leq x'$, $x' \leq x$, $y \leq y'$ and $y' \leq y$. If $x \leq y$, then $x' \leq x \leq y \leq y'$ and this shows that $\leq$ defined on $\overline{P}$ is well-defined. By inheriting the symmetric and transitive properties of $\leq$ on $P$, $\leq$ is also symmetric and transitive on $\overline{P}$. If $\overline{x} \leq \overline{y}$ and $\overline{y} \leq \overline{x}$, then $x \leq y$ and $y \leq x$. This shows that $x \sim y$ i.e., $\overline{x} = \overline{y}$. Hence $\pi$ is order-preserving.

Finally we check the universal property: if $\overline{f} : (\overline{P},\leq) \rightarrow (Q,\preceq)$ is order-preserving, then we define $\overline{f} : \overline{P} \rightarrow Q$ by $\overline{x} \mapsto f(x)$. If $\overline{x} = \overline{y}$, then $x \leq y$ and $y \leq x$. Because $f$ is order-preserving, we obtain $f(x) \preceq f(y)$ and $f(y) \preceq f(x)$. Because  $(Q,\preceq)$ is a poset, we have $f(x) = f(y)$. It is clear that $\overline{f} \circ \pi = f$ and $\overline{f}$ is uniquely determined by $f$.
\end{proof}
Now we can define \textit{Alexandrov topology} on pre-ordered sets $(P,\leq)$.
\begin{Pdef}
Let $(P,\leq)$ be a pre-ordered set. A subset $U$ of $P$ is called \textbf{Alexandrov open} if it satisfies the following property: If $x \in U$ and $x \leq y$, then $y \in U$. All Alexandrov open subsets of $P$ form a topology on $P$, with is called the \textbf{Alexandrov topology} on $(P,\leq)$.
\end{Pdef}
\begin{proof}
Let $\mathfrak{A}$ be the set of all Alexandrov open subsets of $P$, we check that $\mathfrak{A}$ forms a topology on $P$. By definition, $\emptyset, P \in \mathfrak{A}$. For $U, V \in \mathfrak{A}$ and $x, y \in P$ satisfy $x \in U \cap V$ and $x \leq y$, then $y \in U$ and $y \in V$ since $U$ and $V$ are Alexandrov open. Finally, if $\{ U_i \}_{i \in I}$ is a collection of Alexandrov open subsets of $P$ and $x, y \in P$ satisfy $x \in \bigcup_{i \in I} U_i$ and $x \leq y$, then $x \in U_j$ for some $j \in I$. This shows that $y \in U_j \subseteq \bigcup_{i \in I} U_i$, as desired.
\end{proof}
The following proposition shows that the Alexandrov topologies naturally extend the order-preserving functions to be continuous functions between Alexandrov topological spaces.
\begin{prop.}
Let $(P, \leq)$ and $(Q, \preceq)$ be pre-ordered sets. Suppose $f : (P, \leq) \rightarrow (Q, \preceq)$ is a ordering preserving function and $\mathfrak{A}_P$, $\mathfrak{A}_Q$ are Alexandrov topologies on $P$ and $Q$ respectively, then $f : (P, \mathfrak{A}_P) \rightarrow (Q, \mathfrak{A}_Q)$ is continuous.
\end{prop.}
\begin{proof}
Suppose $U$ is an open subset of $Q$. If $x, y \in P$ satisfy $x \in f^{-1}(U)$ and $x \leq y$, then the order-preserving property of $f$ shows that $f(x) \leq f(y)$. Because $f(x) \in U$ and $U$ is Alexandrov open, we deduce that $f(y) \in U$ i.e., $y \in f^{-1}(U)$, hence $f^{-1}(U)$ is Alexandrov open in $P$. 
\end{proof}
Here we give an open basis for a Alexandrov topological space.
\begin{Pdef}
\label{Proposition : Standard open base of Alexandrov topology}
Let $(P,\leq)$ be a pre-ordered set and $\mathfrak{A}$ be the Alexandrov topology on $P$. For every $x \in P$, define the \textbf{open star} at $x$ is the set 
\begin{equation*}
U_x := \{ y \in P : x \leq y \}.
\end{equation*}
Then $U_x \subseteq U_y$ if and only if $y \leq x$ and the family $\mathfrak{B} = \{ U_x : x \in P \}$ is an open basis for $\mathfrak{A}$. 

Similarly, the \textbf{closure} of $x$ is defined by $\overline{x} = \{ y \in P : y \leq x \}$, which is the closure of $\{ x \}$ in $(X,\mathfrak{A})$.
\end{Pdef}
\begin{proof}
It is clear that $U_x \subseteq U_y$ if and only if $y \leq x$. For $x \in P$, $y \in U_x$ and $z \in P$ satisfying $y \leq z$, we have $x \leq y \leq z$ hence $x \leq z$ i.e., $z \in U_x$, hence $U_x$ is Aleandrov open. If $U$ is an Alexandrov open subset of $P$ and $x \in U$, then $U_x \subseteq U$; this shows that $U = \bigcup_{x \in U} U_x$ since $x \in U_x$ whenever $x \in P$. If $x, y \in P$ and $z \in U_x \cap U_y$, then for every $w \in U_z$, we have $x \leq z \leq w$ and $y \leq z \leq w$. Hence $w \in U_x \cap U_y$ and we deduce that $U_z \subseteq U_x \cap U_y$. Therefore, $\mathfrak{B}$ is an open basis for $\mathfrak{A}$. 

Next we show that the closure of $x$ is the Alexandrov closure of the set $\{ x \}$ in $P$ by verifying the following two facts:
\begin{itemize}
\item $\overline{x}$ is Alexandrov closed: For $z \in P \setminus \overline{x}$ and $w \in P$ satisfying $z \leq w$. If $w \notin P \setminus \overline{x}$ i.e. $w \in \overline{x}$, then $w \leq x$. This shows that $z \leq x$ i.e. $z \in \overline{x}$, a contradiction;
\item If $x \in C$ and $C$ is Alexandrov closed, then for every $y \in \overline{x}$, we have $y \leq x$. If $y \in P \setminus C$, then $x \in P \setminus C$ since $P \setminus C$ is Alexandrov open, a contradiction. In other words, $\overline{x}$ must be a subset of $C$. 
\end{itemize}
Therefore, the set $\overline{x} := \{ y \in P : y \leq x \}$ is really the Alexandov closure of $\{x\}$. 
\end{proof}
Now we can view every poset $(P,\leq)$ as a topological space. A natural question is how do we modify a cellular sheaf to a sheaf on the topological space $(P,\mathfrak{A})$. As we proved in Proposition \ref{Proposition : Standard open base of Alexandrov topology}, we know that $U_x \subseteq U_y$ if and only if $y \leq x$. Therefore, if $\shG$ is a cellular sheaf on the poset $(P,\leq)$, the assignment $U_p \mapsto \shG(p)$ is well-defined cause $U_p = U_q$ if and only if $p = q$\footnote{Here we use the anti-symmetry of partial order $\leq$ for $p \leq q$ and $q \leq p$ implying $p = q$. Therefore, the well-definedness of the assignment may not hold if $(P,\leq)$ is just a pre-ordered set.}. By defining $\rho_{U_pU_q} = \rho_{pq}$ for $p \leq q$, this assignment leads commutative diagrams
\begin{equation}
\label{Eq : Diagram of Restriction of maps of Cellular Sheaves-OPEN BASE FORM}
\xymatrix{
\shG(U_p) \ar[r]^{\rho_{U_pU_q}} \ar[dr]_{\rho_{U_pU_r}} & \shG(U_q) \ar[d]^{\rho_{U_qU_r}} \\
& \shG(U_r)
}
\end{equation} 
whenever $p \leq q \leq r$. These are compatible to diagrams of sheaves on topological spaces as in~\eqref{Eq : Diagram of Restriction of maps of Sheaves}. However, it is not sufficient to define objects and morphisms only on $U_p$'s. To achieve this, we introduce the concept of base sheaf defined on a topological space.
\section{Base Sheaves}
\label{sec: Base Sheaves}
Materials in the section are famous exercises in textbooks of algebraic geometry or commutative algebra (\textit{e.g.}~\cite{AM, QingLiu}). A well-known application is using base sheaves to construct affine scheme on the spectrum of a ring~\cite{AM, QingLiu, Ueno, Hartshorne, LeiFu}. In this note, we go through all details of base sheaves as a friendly invitation for beginners in algebraic geometry. 

Let $X$ be a topological space and $\mathfrak{B}$ be an open basis for $X$. A $\mathfrak{B}$\textit{-presheaf} $\shF$ of Abelian groups (resp. rings, algebras, vector spaces, ... etc.) on $X$ consists of the following data :
\begin{itemize}
\itemsep = -1pt
\item For every $U \in \mathfrak{B}$, $\shF$ assigns $U$ to an Abelian group (resp. ring, algebra, vector space, ... etc.) $\shF(U)$ .
\item For every pair $V \subseteq U$ of members in $\mathfrak{B}$, there is a group (resp. ring, algebra, vector space, ... etc.) homomorphism $\rho_{UV} : \shF(U) \rightarrow \shF(V)$, called the \textit{restriction map}. Similarly, we use $s|_V$ to denote $\rho_{UV}(s)$ whenever $s \in \shF(U)$.
\end{itemize}
These data satisfy the following two properties :
\begin{enumerate}
\itemsep = -1pt
\item For any $U \in \mathfrak{B}$, $\rho_{UU} = \id_{\shF(U)}$.
\item For every inclusions $W \subseteq V \subseteq U$ of members in $\mathfrak{B}$, $\rho_{UW} = \rho_{VW} \circ \rho_{UV}$. In other words, for any inclusions $W \subseteq V \subseteq U$, the following diagram is commutative.
\begin{equation}
\xymatrix{
\shF(U) \ar[r]^{\rho_{UV}} \ar[dr]_{\rho_{UW}} & \shF(V) \ar[d]^{\rho_{VW}} \\
& \shF(W)
}
\end{equation}  
\end{enumerate}
A $\mathfrak{B}$-presheaf $\shF$ of Abelian groups (resp. rings, modules, algebras, ... etc.)  on $X$ is a \textit{$\mathfrak{B}$-sheaf} if it satisfies the following additional conditions :
\begin{itemize}
\itemsep = -1pt
\item[$\rm (a)$] Let $U \in \mathfrak{B}$ and $s, t \in \shF(U)$ be sections. If $U$ is covered by $\lbrace U_i \rbrace_{i \in I} \subseteq \mathfrak{B}$ (i.e., $U = \bigcup_{i \in I} U_i$) such that $s|_{U_i} = t|_{U_i}$ for all $i \in I$, then $s = t$.
\item[$\rm (b)$] Suppose $U \in \mathfrak{B}$ and $U$ is covered by $\lbrace U_i \rbrace_{i \in I} \subseteq \mathfrak{B}$ with local sections $s_i \in \shF(U_i)$ such that for every $i, j \in I$, $s_i|_{V} = s_j|_{V}$ whenever $V \in \mathfrak{B}, V \subseteq U_i \cap U_j$. Then there exists a section $s \in \shF(U)$ such that $s|_{U_i} = s_i$ for any $i \in I$. We note that condition $\rm (a)$ implies that the existence of such $s$ is unique. 
\end{itemize}
In other words, a $\mathfrak{B}$-presheaf $\shF$ of Abelian groups is a $\mathfrak{B}$-sheaf if and only if the sequence \begin{equation}
\label{Equation : Exact seq of B-sheaf}
\xymatrix@+1.0em{
				0 \ar[r]_{} & \shF(U)
				\ar[r]^{}
                & \prod_{i \in I} \shF(U_i) 
                \ar[r]^{}
                & \prod_{(i,j) \in I^2} \prod_{V \in \mathfrak{B}, V \subseteq U_i \cap U_j} \shF(V)
}
\end{equation}
equipped with canonical maps is exact for every collection $\{ U_i \}_{i \in I}$ of elements in $\mathfrak{B}$ which covers arbitrary $U$ in $\mathfrak{B}$.

To extend a $\mathfrak{B}$-sheaf to a sheaf on $X$, for arbitrary non-empty open subset $U$ of $X$, we define $I(U) = \{ V : V \in \mathfrak{B}, V \subseteq U \}$. Then $(I(U), \leq)$ is a poset and we define
\begin{equation*}
\shF^+(U) = \mathop{\lim_{\longleftarrow}}_{V \in I(U)} \shF(V) := \left\lbrace (s_V)_{V \in I(U)} \in \prod_{V \in I(U)} \shF(V) : \rho_{VW}(s_V) = s_W \ \forall \ W \subseteq V \right\rbrace 
\end{equation*}
and define $\rho_{UU'}^+ : \shF^+(U) \rightarrow \shF^+(U')$ to be the canonical homomorphism for $U \subseteq U'$ being open in $X$. Then the extended assignment $\shF$ is a presheaf on $X$.
\begin{remark}

If $V \in \mathfrak{B}$, then $\shF^+(V) \simeq \shF(V)$ canonically since each $I(V)$-tuple $(s_W)_{W \in I(V)}$ is unique determined by $s_V$. On the other hand, suppose $W \subseteq V$ are members in $\mathfrak{B}$ and $\mathbf{s} = (s_O)_{O \in I(V)}$, then the $W^{\rm th}$ component of $\rho_{VW}^+(\mathbf{s})$ is the $W^{\rm th}$ component of $\mathbf{s}$. By the constrain of projective limit, which is $\rho_{VW}(s_V)$. By the identification $\shF^+(V) \simeq \shF(V)$, $\rho_{VW}^+$ and $\rho_{VW}$ are coincide. 
\end{remark}
Alternatively, the $\shF^+(U)$ can be viewed as the set of all functions $\mathbf{s}$ from $I(U)$ to $\coprod_{V \in I(X)} \shF(V)$, which satisfy the following two properties:
\begin{itemize}
\item[$(\rm i)$] For every $V \in I(U)$, $\mathbf{s}(V) \in \shF(V)$;
\item[$(\rm ii)$] For every $V, W \in I(U)$ satisfying $W \subseteq V$, $\rho_{VW}(\mathbf{s}(V)) = \mathbf{s}(W)$.
\end{itemize}
Up to this setting, the canonical homomorphism $\rho_{UU'}^+ : \shF^+(U) \rightarrow \shF^+(U')$ is just the restriction of functions $\mathbf{s} \mapsto \mathbf{s}|_{I(U')}$.
\begin{lemma}
\label{Lemma : I(U) property}
Let $X$ be a topological space, $\mathfrak{B}$ an open basis for $X$, and $I(U) = \{ V \in \mathfrak{B} : V \subseteq U \}$. Then the following formulas hold:
\begin{itemize}
\item[$(\rm i)$] $U_1 \subseteq U_2$ if and only if $I(U_1) \subseteq I(U_2)$;
\item[$(\rm ii)$] $U_1 = U_2$ if and only if $I(U_1) = I(U_2)$;
\item[$(\rm iii)$] $I(U_1 \cap U_2) = I(U_1) \cap I(U_2)$.
\end{itemize}
\end{lemma}
\begin{proof}
(i) It is clear that $U_1 \subseteq U_2$ implies $I(U_1) \subseteq I(U_2)$. Conversely, if $I(U_1) \subseteq I(U_2)$, then $U_1 = \bigcup_{j \in J} V_{j}$ where $V_{j} \in \mathfrak{B}$ for each $j \in J$ since $\mathfrak{B}$ is an open basis for $X$. Because $V_j \in I(U_1) \subseteq I(U_2)$ for each $j$, we obtain $U_1 \subseteq U_2$. (ii) follows from (i) immediately.

(iii) By (i), $I(U_1 \cap U_2) \subseteq I(U_1) \cap I(U_2)$ because $I(U_1 \cap U_2) \subseteq I(U_1)$ and $I(U_1 \cap U_2) \subseteq I(U_2)$. Conversely, if $V \in I(U_1) \cap I(U_2)$, then $V \in \mathfrak{B}$ and $V \subseteq U_1, V \subseteq U_2$. This shows that $V \in \mathfrak{B}$ and $V \subseteq U_1 \cap U_2$ i.e., $V \in I(U_1 \cap U_2)$.
\end{proof}
In fact, $\shF^+$ is a sheaf, it can be verified as follow:
\begin{itemize}
\item[$\rm (a)$] Let $\{ U_i \}_{i \in I}$ be an open cover of an open subset $U \subseteq X$ and $\mathbf{s} \in \shF^+(U)$ satisfy $\mathbf{s}|_{I(U_i)} = 0$ for every $i \in I$. For $V \in I(U)$, each $V \cap U_i$ is open in $X$, hence it can be covered by elements in $\mathfrak{B}$. Hence we can assume $V = \bigcup_{k \in K} W_k$ where each $W_k$ is an element in $\mathfrak{B}$ and $W_k$ is contained in some $U_{i(k)}$ with $i(k) \in I$. For each $k \in K$, $\mathbf{s}|_{I(W_k)} = \mathbf{s}|_{I(U_{i(k)})}|_{I(W_{k})} = 0$, this shows that $\rho_{VW_k}(\mathbf{s}(V)) = \mathbf{s}(W_k) = \mathbf{s}|_{I(W_k)}(W_k) = 0$. By (a) of $\mathfrak{B}$-sheaf of $\shF$, we deduce that $\mathbf{s}(V) = 0$. Hence $\mathbf{s} = 0$.
\item[$\rm (b)$] Let $\{ U_i \}_{i \in I}$ be an open cover of an open subset $U \subseteq X$ with local sections $\mathbf{s}_i \in \shF^+(U_i)$ satisfy $\mathbf{s}_i|_{I(U_i) \cap I(U_j)} = \mathbf{s}_j|_{I(U_i) \cap I(U_j)}$ whenever $i, j \in I$. Let $V$ be an element in $I(U)$, then
\begin{equation}
\label{Equation : Totally cover}
V = \bigcup_{i \in I} V \cap U_i = \bigcup_{i \in I} \left( \bigcup_{W \in I(V \cap U_i)} W \right)
\end{equation}
For each $W$ in \eqref{Equation : Totally cover}, $W \subseteq U_{i(W)}$ for some $i(W) \in I$, define
\begin{equation*}
t_W = \mathbf{s}_{i(W)}(W) \in \shF(W).
\end{equation*}
Note that $t_W$ is independent to the choice of $i(W)$. That is, if $W \subseteq U_j$ with $j \in I$, then $W \in I(U_j) \cap I(U_{i(W)})$ and hence $\mathbf{s}_j(W) = \mathbf{s}_{i(W)}(W) = t_W$ by compatibility of local sections. Then for $W_1, W_2$ in \eqref{Equation : Totally cover}, if $W \in I(W_1 \cap W_2) \subseteq I(U_{i(W_1)} \cap U_{i(W_2)}) = I(U_{i(W_1)}) \cap I(U_{i(W_2)})$ (by Lemma \ref{Lemma : I(U) property}), we have
\begin{equation*}
\begin{split}
\rho_{W_1W}(t_{W_1}) &= {\mathbf{s}}_{i(W_1)}(W) = {\mathbf{s}}_{i(W_2)}(W) = \rho_{W_2W}(t_{W_2}). 
\end{split}
\end{equation*}
By (b) of $\mathfrak{B}$-sheaf structure of $\shF$, there is a $t \in \shF(V)$,  such that 
\begin{equation}
\label{Equation : Key Global}
\rho_{VW}(t) = t_W = \mathbf{s}_{i(W)}(W)
\end{equation}
for all $W$ in \eqref{Equation : Totally cover}. Suppose $V, V' \in I(U)$ satisfy $V' \subseteq V$ and $t \in \shF(V), t' \in \shF(V')$ are derived by previous arguments respectively. If $W \in I(V') \subseteq I(V)$ and $W \subseteq U_{i(W)}$ for some $i(W) \in I$, we have
\begin{equation*}
\begin{split}
\rho_{V'W}(t') =\mathbf{s}_{i(W)}(W) = \rho_{VW}(t) = \rho_{V'W}(\rho_{VV'}(t))
\end{split}
\end{equation*}
via equation \eqref{Equation : Key Global}. By (a) of $\mathfrak{B}$-sheaf structure of $\shF$, $\rho_{VV'}(t) = t'$. This shows that all $(t,V)$ really define an element $\mathbf{s}$ in $\shF^+(U)$ such that $\mathbf{s}|_{U_i} = \mathbf{s}_i$. 
\end{itemize}
Now we can successfully extend a $\mathfrak{B}$-sheaf as a sheaf on the whole space. The following proposition tells us that, this extension would not change the local informations.
\begin{prop.}
\label{Prop. B-sheaf stalks at P equals the stalk of the induced sheaf at P}
Let $(X,\mathfrak{T})$ be a topological space, $\mathfrak{B}$ an open basis for $X$, and $\shF$ and $\mathfrak{B}$-sheaf of Abelian groups (resp. rings, modules, ... etc.) . Let $\shF^+$ be the sheaf on $X$ which is induced by $\shF$, then $\shF^+_P \simeq \shF_P$ as groups (resp. rings, modules, ... etc.)  for every $P \in X$.
\end{prop.}
\begin{proof}
As we remarked, we view $\shF(U) = \shF^+(U)$ whenever $U \in \mathfrak{B}$. Because $\mathfrak{B}$ is a subset of $\mathfrak{T}$ such that $(\mathfrak{B}, \supseteq)$ is also a direct set, for each $P \in X$, we have canonical group (resp. ring) homomorphism
\begin{equation*}
\phi : \mathop{\lim_{\longrightarrow}}_{U \in  \mathfrak{B}, P \in U} \shF(U) \longrightarrow \mathop{\lim_{\longrightarrow}}_{U \in  \mathfrak{T}, P \in U} \shF(U) = \shF^+_P
\end{equation*}
Because each element in $\shF_P^+$ can be represented by $\mathbf{s}_P$ for some $\mathbf{s} \in \shF^+(U)$ with $P \in U$ and $U \in \mathfrak{B}$, and $\mathbf{s}$ corresponds to $\mathbf{s}(U) \in \shF(U)$, $\phi$ is bijective. 
\end{proof}
\section{Universal Property of Base Sheaves}
\label{sec : Universal Property of Base Sheaves}
As we learned in linear algebra, if $V$ and $W$ are vector spaces over a field and $\beta$ is a basis for $V$, to define a linear transformation $T : V \rightarrow W$, it is sufficient to define a mapping $\beta \rightarrow W$ and linearly extend this mapping. An analogue version of extending morphisms on base sheaves is listed as follow, which is called \textit{the universal property of base sheaves}.
\begin{prop.}
[\cite{QingLiu} Exercise 2.7, the first part]
\label{Proposition : Universal property of B-sheaf}
Let $\mathfrak{B}$ be an open basis on a topological space $(X, \mathfrak{T})$. Let $\shF$ and $\shG$ be sheaf and $\mathfrak{B}$-sheaf of Abelian groups (resp. rings)  on $X$ respectively. Suppose that for every $V \in \mathfrak{B}$ there exists a homomorphism $\alpha(V) : \shF(V) \rightarrow \shG(V)$ which is compatible with restrictions. Then it extends in a unique way to a homomorphism of sheaves $\alpha^+ : \shF \rightarrow \shG^+$. 
\end{prop.}
\begin{proof}
For every $U$ open in $X$ and $V, W \in I(U)$ with $W \subseteq V$, we have the following commutative diagram:
\begin{equation*}
\xymatrix@+1.5em{
& \shF(U) \ar[ldd]_{\alpha(V) \circ \rho_{UV}}\ar[rdd]^{\alpha(W) \circ \rho_{UW}}\ar@{-->}[d]^{} & \\
& \shG^+(U) \ar[ld]^{\pi_{UV}} \ar[rd]_{\pi_{UW}} & \\
\shG(V) \ar[rr]_{\sigma_{VW}} & & \shG(W) }
\end{equation*}
where $\rho_{\bullet\bullet}$ and $\sigma_{\bullet\bullet}$ denote the restriction maps of $\shF$ and $\shG$ respectively. Note that it is commutative because $\sigma_{VW} \circ \alpha(V) \circ \rho_{UV} = \alpha(W) \circ \rho_{VW} \circ \rho_{UV} = \alpha(W) \circ \rho_{UW}$. By the universal property of projective limits, there is a unique $\alpha^+(U) : \shF(U) \rightarrow \shG^+(U)$ such that $\pi_{UW} \circ \alpha^+(U) = \alpha(W) \circ \rho_{UW}$ whenever $W \in I(U)$. Note that the existence of $\alpha^+$ is a morphism of sheaves because for arbitrary $U' \subseteq U$ and $W \in I(U')$,
\begin{equation*}
\begin{split}
\pi_{U'W} \circ \sigma_{UU'}^+ \circ \alpha^{+}(U) &= \pi_{UW} \circ \alpha^+(U) = \alpha(W) \circ \rho_{UW}, \\
\pi_{U'W} \circ \alpha^{+}(U') \circ \rho_{UU'} &= \alpha(W) \circ \rho_{U'W} \circ \rho_{UU'} = \alpha(W) \circ \rho_{UW}.
\end{split}
\end{equation*}
Because elements in $\shG^+(U')$ are uniquely determined by projections $\pi_{U'
W}$, $W \in I(U')$, we obtain $\sigma_{UU'}^+ \circ \alpha^{+}(U) = \alpha^{+}(U') \circ \rho_{UU'}$. Finally, when $U \in \mathfrak{B}$, $\pi_{UU} \circ \alpha^+(U) = \alpha(U) \circ \rho_{UU} = \pi_{UU} \circ \alpha^+(U) = \alpha(U)$. In this case, $\pi_{UU} : \shG^+(U) \rightarrow \shG(U)$ is the canonical isomorphism and hence $\alpha^+(U)$ and $\alpha(U)$ are coincide.

For uniqueness, because $\beta(U) : \shF(U) \rightarrow \shG^+(U)$ satisfies $\beta(U') \circ \rho_{UU'} = \sigma_{UU'}^+ \circ \beta(U)$ whenever $U' \subseteq U$ implies that $\beta(W) \circ \rho_{UW} = \sigma_{UW}^+ \circ \beta(U)$ whenever $W \in I(U)$. Hence
\begin{equation*}
\sigma_{UW}^+ \circ \alpha(U) = \alpha(W) \circ \rho_{UW} = \beta(W) \circ \rho_{UW} = \sigma_{UW}^+ \circ \beta(U.)
\end{equation*}
Because each element in $\shG^{+}(U)$ is uniquely determined by values sent by $\rho_{UW}^+$ with $W \in I(U)$, we conclude $\beta(U) = \alpha(U)$.
\end{proof}
\begin{coro}
Let $\mathfrak{B}$ be an open basis on a topological space $(X, \mathfrak{T})$. Let $\shF$ be a sheaf of Abelian groups (resp. rings, modules, ... etc.) on $X$. Then the sheaf $\shF^+$  induced by its $\mathfrak{B}$-sheaf structure is isomorphic to $\shF$ canonically. 
\end{coro}
\begin{proof}
By Proposition \ref{Proposition : Universal property of B-sheaf}, there is a unique morphism $\beta : \shF \rightarrow \shF^+$ where $\beta^+(U) : \shF(U) \rightarrow \shF^+(U)$ is defined by $s \mapsto (s|_V)_{V \in I(U)}$. Clearly, it is an isomorphism by the construction of sheaves induced by $\mathfrak{B}$-sheaves.
\end{proof}
\begin{prop.}
[\cite{QingLiu} Exercise 2.7, the second part]
Let $\mathfrak{B}$ be an open basis on a topological space $(X, \mathfrak{T})$. Let $\shF, \shG$ be two sheaves on $X$. Suppose that for every $U \in \mathfrak{B}$ there exists a homomorphism $\alpha(U) : \shF(U) \rightarrow \shG(U)$ which is compatible with restrictions. Then it extends in a unique way to a homomorphism of sheaves $\alpha : \shF \rightarrow \shG$. Moreover, if $\alpha$ is surjective (resp. injective) for every $U \in \mathfrak{B}$, then $\alpha$ is surjective (resp. injective).
\end{prop.}
\begin{proof}
By the universal property of $\mathfrak{B}$-sheaves, there is a unique morphism $\alpha^+ : \shF \rightarrow \shG^+$ which extend $\alpha$. Because $\shG^+ \simeq \shG$ canonically, the uniqueness and existence follow.

Let $P$ be a point in $X$. If $\alpha(U) : \shF(U) \rightarrow \shG(U)$ is surjective for every $U \in \mathfrak{B}$, then for each $t_P \in \shG_P$ pick an $U \in \mathfrak{B}$ such that $P \in U$. Because $\alpha(U)$ is onto, there is an $s \in \shF(U)$ such that $\alpha(U)(s) = t$ and we have $\alpha_P(s_P) = \alpha(U)(s)_P = t_P$. On the other hand, if $\alpha(U) : \shF(U) \rightarrow \shG(U)$ is one-to-one and $\alpha_P(s_P) = 0$ with $s \in \shF(U)$ for some $U \in \mathfrak{B}$, then there an $V \in \mathfrak{B}$, $V \subseteq U$ such that $\alpha(V)(s|_V) = \alpha(U)(s)|_V = 0$. Because $\alpha(V)$ is one-to-one, $s|_V = 0$ and we conclude that $s_P = 0$.
\end{proof}
\section{Cellular Sheaves on Alexandrov Spaces}
Now we can show how do we extend a cellular sheaf $\shG$ of $A$-modules on a poset $(P, \leq)$ to a sheaf $\shG^+$ of $A$-modules on the Alexandorv topological space $(P, \mathfrak{A})$. By the $\mathfrak{B}$-sheaf structure, it is sufficient to show that the assignment (still denoted by $\shG$)
\begin{equation}
U_p \longmapsto \shG(U_p) := \shG(p) \ \ {and} \ \ \rho_{U_pU_q} := \rho_{pq}   
\end{equation}
whenever $p \leq q$ in $P$ defines a $\mathfrak{B}$-sheaf of $A$-modules on $(P, \mathfrak{A})$. Because $\shG$ is a cellular sheaf on $(P, \leq)$, the pre-sheaf structure holds naturally. Here we check $\shG$ satisfies properties (a) and (b) as axioms of base sheaves.
\begin{itemize}
    \item[\rm (a)] Suppose $p \in P$, $s \in \shG(U_p)$ and $\{ U_{p_i} : p_i \in P, i \in I \}$ is an open cover of $U_p$ such that $s|_{U_{p_i}} = 0$ for each $i \in I$. Because $U_{p_i} \subseteq U_p$, we have $p \leq p_i$ for all $i \in I$. On the other hand, there is a $j \in I$ such that $p \in U_{p_j}$, this shows that $p_j \leq p \leq p_j$. In particular, $s = s|_{U_p} = s|_{U_{p_j}} = 0$;
    \item [\rm (b)] Suppose $p \in P$ and $\{ U_{p_i} : p_i \in P, i \in I \}$ is an open cover of $U_p$. As we proved in (a), $p = p_k$ for some $k \in I$. If $s_i \in \shG(U_{p_i})$ are local sections such that $s_i|_{U_{p_i} \cap U_{p_j}} = s_j|_{U_{p_i} \cap U_{p_j}}$ whenever $i, j \in I$, then $s_k \in \shG(U_{p_k}) = \shG(U_{p})$ and $s_k|_{U_{p_i}} = s_k|_{U_{p} \cap U_{p_i}} = s_k|_{U_{p_k} \cap U_{p_i}} = s_i|_{U_{p_k} \cap U_{p_i}} = s_i|_{U_{p} \cap U_{p_i}} = s_i|_{U_{p_i}} = s_i$.
\end{itemize}
There, $\shG$ is a $\mathfrak{B}$-sheaf of $A$-modules on $(P, \mathfrak{A})$ where $\mathfrak{B} = \{ U_p : p \in P \}$ is an open basis for $\mathfrak{A}$. By the discussions in Section~\ref{sec: Base Sheaves} and~\ref{sec : Universal Property of Base Sheaves}, $\shG$ uniquely determines a sheaf $\shG^+$ of $A$-modules on the topological space. 
\begin{theorem}
[\cite{justinCurryPHDThesis}~Theorem 4.2.10]
\label{Theorem : Stalks of a cellular sheaf are really stalks}
Let $\shG$ be a cellular sheaf of $A$-modules on a poset $(P,\leq)$ and $\shG^+$ be the sheaf on $(P, \mathfrak{A})$ induced by $\shG$. Then for every $p \in P$, the stalk $\shG^+_p$ of $\shG^+$ at $p$ is canonically isomorphic to $ \shG(p) = \shG(U_p)$. 
\end{theorem}
\begin{proof}
By Proposition \ref{Prop. B-sheaf stalks at P equals the stalk of the induced sheaf at P}, $\shG^+_p \simeq \shG_p$ canonically. Define $\phi : \shG(U_p) \rightarrow \shG_p$ by $s \mapsto s_p$, where $s_p$ denotes the stalk of local section $s \in \shF(U_p)$ at $p$. It is clear that $\phi$ is an $A$-module homomorphism and $\phi$ is onto because every open neighbourhood of $p$ in $(P,\mathfrak{A})$ must contains $U_p$. Finally, if $s \in \ker(\phi)$ with $s \in \shF(U_p)$, then $s_p = 0$ implies that $s|_{U_q} = 0$ where $U_q \subseteq U_p$ is an open neighbourhood of $p$ by the definition of direct limit. However, $p \in U_q$ implies that $U_p \subseteq U_q \subseteq U_p$, this forces that $s = 0$. 
\end{proof}
In previous works~\cite{RobinsonNyquist, RobinsonHunting, robinson2014topological, ghrist2019}, authors usually directly defined $\shG(p)$ as the ``stalk'' of the cellular sheaf $\shG$ at $p$. The theorem gives an explanation that it is compatible to the definition of stalks for the induced sheaf on the Alexandrov topological space (cf.~\cite{justinCurryPHDThesis}~Remark 4.2.11).  
\section{Examples}
In the end of this note, we exhibit some basic examples of cellular sheaves. For convenience, we usually use diagram and arrows to represent posets. That is, $p \leq q$, then we can represent it by $p \rightarrow q$. For instance, the poset $(\bbZ, \leq)$ can be represented by
\begin{equation}
\label{Eq : Diagrams of Integers}
    \cdots \rightarrow -2 \rightarrow -1 \rightarrow 0 \rightarrow 1 \rightarrow 2 \rightarrow \cdots.
\end{equation}
Note that the self loops for representing $x \leq x$ are omitted in the diagrams. This notations were used frequently in the category theory~\cite{justinCurryPHDThesis, CateMacLane, CateAwodey, CateRiehl}. In the example \eqref{Eq : Diagrams of Integers}, a cellular sheaf $\shG$ of $A$-modules would be a sequence
\begin{equation}
\label{Eq : Diagrams of Cellular sheaves}
    \cdots \rightarrow \shG(-2) \rightarrow \shG(-1) \rightarrow \shG(0) \rightarrow \shG(1) \rightarrow \shG(2) \rightarrow \cdots.
\end{equation}
of $A$-modules and $A$-module homomorphisms. 
\begin{exam.}
Let $(P, \leq)$ and $\shG$ be poset and cellular sheaf of $A$-modules respectively represented by
\begin{equation}
\xymatrix@+2.0em{
				& q_1
				\ar[r]_{}
                & r 
                \\
        		& p \ar[u]^{} \ar[r]_{}
				& q_2 \ar[u]^{}
} 
\xymatrix@+2.0em{
				& \shG(q_1)
				\ar[r]^{\rho_{q_1r}}
                & \shG(r) 
                \\
        		& \shG(p) \ar[u]^{\rho_{pq_1}} \ar[r]_{\rho_{pq_2}}
				& \shG(q_2) \ar[u]_{\rho_{q_2r}}
}
\end{equation}
Then the local section of $U := U_{q_1} \cup U_{q_2}$ is the $A$-module $\{ (s_1, s_2) \in \shG(q_1) \times \shG(q_2) : s_1|_{U_r} = s_2|_{U_r} \}$.
\end{exam.}
\begin{exam.}
Let $(P, \leq)$ and $\shG$ be poset and cellular sheaf of $A$-modules respectively represented by
\begin{equation}
\xymatrix@+1.5em{
				& p_1
				\ar[dr]_{}
                &  
                \\
        		& p_k \ar@{.}[u]^{} \ar[r]_{}
				& q 
				\\
        		& p_n \ar@{.}[u]^{} \ar[ur]_{}
				& 
} 
\xymatrix@+1.5em{
				& \shG(p_1)
				\ar[dr]^{\rho_{p_1q}}
                &  
                \\
        		& \shG(p_k) \ar@{.}[u]^{} \ar[r]_{\rho_{p_kq}}
				& \shG(q) 
				\\
        		& \shG(p_n) \ar@{.}[u]^{} \ar[ur]_{\rho_{p_nq}}
				& 
}
\end{equation}
Then the local section of $U := U_{p_1} \cup \cdots \cup U_{p_n}$ is the $A$-module $\{ (s_i)_{i = 1}^n \in \prod_{i = 1}^n \shG(p_i) : s_i|_{U_q} = s_j|_{U_q} \ \forall i,j \}$.
\end{exam.}
\begin{exam.}
Let $(P, \leq)$ and $\shG$ be poset and cellular sheaf of $A$-modules respectively represented by
\begin{equation}
\xymatrix@+1.5em{
                & 
				& p_1 \ar[dl]_{}
				\ar[dr]_{}
                &  
                \\
                & q_1   
        		& p_k \ar@{.}[u]^{} \ar[r]_{} \ar[l]_{}
				& q_2 
				\\
				&
        		& p_n \ar@{.}[u]^{} \ar[ur]_{} \ar[ul]_{}
				& 
} 
\xymatrix@+1.5em{
                & 
				& \shG(p_1) \ar[dl]_{\rho_{p_1q_1}}
				\ar[dr]^{\rho_{p_1q_2}}
                &  
                \\
                & \shG(q_1) 
        		& \shG(p_k)  \ar[r]^{\rho_{p_kq_2}} \ar@{.}[u]^{} \ar[l]_{\rho_{p_kq_1}}
				& \shG(q_2) 
				\\
				&
        		& \shG(p_n) \ar@{.}[u]^{} \ar[ur]_{\rho_{p_nq_2}}  \ar[ul]^{\rho_{p_nq_1}}
				& 
}
\end{equation}
Then the local section of $U := U_{p_1} \cup \cdots \cup U_{p_n}$ is the $A$-module 
\begin{equation}
    \left \{ (s_i)_{i = 1}^n \in \prod_{i = 1}^n \shG(p_i) : s_i|_{U_{q_k}} = s_j|_{U_{q_k}} \ \forall i,j = 1,2..., n \ \forall k = 1,2 \right \}.
\end{equation}
\end{exam.}
\begin{exam.}
Let $(P, \leq)$ and $\shG$ be poset and cellular sheaf of $A$-modules respectively represented by
\begin{equation}
\xymatrix@+1.5em{
                & 
				& p_1 
				\\
                & q \ar[r]_{} \ar[ur]_{} \ar[dr]_{} 
        		& p_k \ar@{.}[u]^{} 
				\\
				&
        		& p_n \ar@{.}[u]^{} 
} 
\xymatrix@+1.5em{
                & 
				& \shG(p_1) 
				\\
                & \shG(q) \ar[r]^{\rho_{qp_k}} \ar[ur]^{\rho_{qp_1}} \ar[dr]_{\rho_{qp_n}} 
        		& \shG(p_k) \ar@{.}[u]^{} 
				\\
				&
        		& \shG(p_n) \ar@{.}[u]^{} 
} 
\end{equation}
Then the local section of $U := U_{p_1} \cup \cdots \cup U_{p_n}$ and $U_q$ are the $A$-modules $\prod_{i = 1}^n \shG(p_i)$ and $\shG(q)$ respectively.
\end{exam.}
\section{Conclusion}
In this note, we briefly introduced the connection between cellular sheaves defined on posets and induced sheaf structures on the posets equipped with the Alexandrov topology via base sheaf structure. Instead of using Kan extensions, base sheaves provide an elementary approach to explain the connection. On the other hand, because of the structure of Alexandrov topology, any functor from a poset to the category of modules satisfies the conditions (a) and (b) of base sheaves automatically. This topological structures provide a simple way to understand local sections defined on posets, and may have potential in analyzing topological information behind digital data.    
\section*{Acknowledgments}
The author and the student algebraic geometry seminar are supported by the projects MOST 108-2119-M-002-031 and MOST 108-2115-M-003-005-MY2 hosted by the Ministry of Science and Technology in Taiwan.
\bibliographystyle{unsrt}  
\bibliography{references}

\begin{thebibliography}{10}

\bibitem{Hartshorne}
Robin Hartshorne.
\newblock {\em Algebraic Geometry}.
\newblock Springer-Verlag, New York, 1977.

\bibitem{Ueno}
Kenji Ueno.
\newblock {\em Algebraic Geometry 1: From Algebraic Vaieties to Schemes}.
\newblock Translations of Mathematical Monographs, Volume 185, AMS, 1991.

\bibitem{LeiFu}
Lei Fu.
\newblock {\em Algebraic Geometry}.
\newblock Springer and Tsing Hua University Press, Mathematics Series for
  Graduate Students, 2006.

\bibitem{QingLiu}
Qing Liu.
\newblock {\em Algebraic Geometry and Arithmetic Curves}.
\newblock Oxford Graduate Texts in Mathematics, 2006.

\bibitem{bredon1997sheaf}
G.E. Bredon.
\newblock {\em Sheaf Theory}.
\newblock Graduate Texts in Mathematics. Springer New York, 1997.

\bibitem{ghrist2014elementary}
R.~Ghrist.
\newblock {\em Elementary Applied Topology}.
\newblock CreateSpace Independent Publishing Platform, 2014.

\bibitem{robinson2014topological}
M.~Robinson.
\newblock {\em Topological Signal Processing}.
\newblock Mathematical Engineering. Springer Berlin Heidelberg, 2014.

\bibitem{RobinsonHunting}
Michael Robinson.
\newblock Hunting for foxes with sheaves.
\newblock {\em Notices of the American Mathematical Society}, page 661, 05
  2019.

\bibitem{Curry2015}
Justin~Michael Curry.
\newblock Topological data analysis and cosheaves.
\newblock {\em Japan Journal of Industrial and Applied Mathematics}, 2015.

\bibitem{justinCurryPHDThesis}
Justin~Michael Curry.
\newblock {\em Sheaves, Cosheaves and Applications}.
\newblock PhD Thesis, 2014.

\bibitem{Ladkani_2008}
Sefi Ladkani.
\newblock On derived equivalences of categories of sheaves over finite posets.
\newblock {\em Journal of Pure and Applied Algebra}, 212(2):435–451, Feb
  2008.

\bibitem{ghrist2020cellular}
Robert Ghrist and Hans Riess.
\newblock Cellular sheaves of lattices and the tarski laplacian, 2020.

\bibitem{RobinsonNyquist}
Michael Robinson.
\newblock The nyquist theorem for cellular sheaves.
\newblock {\em Sampling Theory and Applications 2013, Bremen, Germany}, 2013.

\bibitem{AM}
M.F. Atiyah and I.G. Macdonald.
\newblock {\em Introduction to Commutative Algebra}.
\newblock Addison-Wesley Publishing Company, 1969.

\bibitem{RonHeld1}
Ron Held.
\newblock {\em Sheaves Over Posets (draft)}.
\newblock Academia Open Resource, 2016.

\bibitem{Ale37}
P.~Alexandroff.
\newblock {\em Diskrete räume}.
\newblock Mat. Sb, 2:501-518, 1937.

\bibitem{Ale47}
P.~S. Alexandrov.
\newblock {\em Combinatorial Topology}.
\newblock Dover Publications, Inc., 31 East 2nd Street, Mineola, NY, 1947.

\bibitem{munkres2018elements}
J.R. Munkres.
\newblock {\em Elements of Algebraic Topology}.
\newblock CRC Press, 2018.

\bibitem{ghrist2019}
Jakob Hansen and Robert Ghrist.
\newblock Toward a spectral theory of cellular sheaves.
\newblock {\em Journal of Applied and Computational Topology}, 2019.

\bibitem{CateMacLane}
Saunders~Mac Lane.
\newblock {\em Categories for the Working Mathematician}.
\newblock Graduate Texts in Mathematics book series (GTM, volume 5),
  Springer-Verlag New York, 1971.

\bibitem{CateAwodey}
Steve Awodey.
\newblock {\em Category Theory}.
\newblock Oxford University Press, 2nd edition, 2010.

\bibitem{CateRiehl}
Emily Riehl.
\newblock {\em Category Theory in Context}.
\newblock Cambridge University Press, 2014.

\end{thebibliography}

\end{document}